\documentclass{agtart_a}
\pdfoutput=1

\usepackage[all]{xy}
\SelectTips{cm}{12}

\title{Dehn surgery, homology and hyperbolic volume}

\author{Ian Agol}
\givenname{Ian}
\surname{Agol}
\address{Department of Mathematics, 
Statistics, and Computer Science (M/C 249)\\\newline
University of Illinois at Chicago\\
851 S Morgan St\\
Chicago, IL 60607-7045\\USA}
\email{agol@math.uic.edu}
\urladdr{}

\author{Marc Culler}
\givenname{Marc}
\surname{Culler}
\email{culler@math.uic.edu}
\urladdr{}

\author{Peter B Shalen}
\givenname{Peter B}
\surname{Shalen}
\email{shalen@math.uic.edu}
\urladdr{}

\volumenumber{6}
\issuenumber{}
\publicationyear{2006}
\papernumber{80}
\startpage{2297}
\endpage{2312}

\doi{}
\MR{}
\Zbl{}

\keyword{hyperbolic manifold}
\keyword{volume}
\keyword{homology}
\keyword{drilling}
\keyword{Dehn surgery}
\subject{primary}{msc2000}{57M50}
\subject{secondary}{msc2000}{57M27}

\received{14 July 2006}
\revised{}
\accepted{1 November 2006}
\published{8 December 2006}
\publishedonline{8 December 2006}
\proposed{}
\seconded{}
\corresponding{}
\editor{CPR}
\version{}

\arxivreference{math.GT/0508208}

\let\xysavmatrix\xymatrix
\def\xymatrix{\disablesubscriptcorrection\xysavmatrix}
\AtBeginDocument{\let\bar\wbar\let\tilde\wtilde\let\hat\what}

\makeatletter
\def\cnewtheorem#1[#2]#3{\newtheorem{#1}{#3}[section]
\expandafter\let\csname c@#1\endcsname\c@para}
\makeatother

\theoremstyle{remark}
\newtheorem{para}{}[section]
\cnewtheorem{remark}[para]{Remark}
\cnewtheorem{remarks}[para]{Remarks}

\theoremstyle{definition}
\cnewtheorem{definition}[para]{Definition}
\cnewtheorem{definitions}[para]{Definitions}
\cnewtheorem{notation}[para]{Notation}
\cnewtheorem{convention}[para]{Convention}

\theoremstyle{plain}
\cnewtheorem{theorem}[para]{Theorem}
\cnewtheorem{lemma}[para]{Lemma}
\cnewtheorem{proposition}[para]{Proposition}
\cnewtheorem{corollary}[para]{Corollary}

\numberwithin{equation}{para}

\newcommand\Number{\begin{para}}
\newcommand\EndNumber{\end{para}}
\newcommand\Definition{\begin{definition}}
\newcommand\EndDefinition{\end{definition}}
\newcommand\Definitions{\begin{definitions}}
\newcommand\EndDefinitions{\end{definitions}}
\newcommand\Theorem{\begin{theorem}}
\newcommand\EndTheorem{\end{theorem}}
\newcommand\Remark{\begin{remark}}
\newcommand\EndRemark{\end{remark}}
\newcommand\Remarks{\begin{remarks}}
\newcommand\EndRemarks{\end{remarks}}
\newcommand\Convention{\begin{convention}}
\newcommand\EndConvention{\end{convention}}
\newcommand\Notation{\begin{notation}}
\newcommand\EndNotation{\end{notation}}
\newcommand\Lemma{\begin{lemma}}
\newcommand\EndLemma{\end{lemma}}
\newcommand\Proposition{\begin{proposition}}
\newcommand\EndProposition{\end{proposition}}
\newcommand\Corollary{\begin{corollary}}
\newcommand\EndCorollary{\end{corollary}}
\newcommand\Proof{\begin{proof}}
\newcommand\EndProof{\end{proof}}
\newcommand\Equation{\begin{equation}}
\newcommand\EndEquation{\end{equation}}
\newcommand\NoProof{}

\newdimen\partindent
\global\partindent=30pt
\newcommand\Parts{\begingroup}
\newcommand\Part[1]{\par\noindent\hangindent\partindent
     \hbox to \partindent{\hskip .5\partindent minus .5\partindent
     {\rm #1}\qua\hfill}\ignorespaces}
\newcommand\Subpart[1]{\par\noindent\hangindent1.5\partindent
     \hbox to 1.5\partindent{\hskip \partindent minus \partindent
     {\rm #1}\qua\hfill}\ignorespaces}
\newcommand\EndParts{\par\endgroup}

\renewcommand\Re{\mathop{\rm Re}}

\newcommand\ee{{\mathcal E}}
\newcommand\calp{{\mathcal Q}}
\newcommand\F{{\Phi}}
\newcommand\eff{f}
\newcommand\omigosh{w}
\newcommand\gimpty{GMT\ }
\newcommand\point{P}
\newcommand\ek{\ee_k}
\newcommand\volthree{{\rm Vol}3}
\newcommand\Isom{{\rm Isom}}
\newcommand\ZZ{{\mathbb Z}}
\newcommand\CC{{\mathbb C}}
\newcommand\QQ{{\mathbb Q}}
\newcommand\HH{{\mathbb H}}
\newcommand\calh{{\mathcal H}}
\newcommand\cals{{\mathcal S}}
\newcommand\dist{\mathop{\rm dist}}
\newcommand\vol{\mathop{\rm vol}}
\newcommand\arcsinh{\mathop{\rm arcsinh}}

\newcommand\tuberad{\mathop{{\rm tuberad}}}
\newcommand\tube{\mathop{{\rm tube}}}
\newcommand\drillc{\mathop{{\rm drill}_C}}
\newcommand\length{\mathop{{\rm length}}}
\newcommand\image{\mathop{{\rm Im}}}
\newcommand\thebound{1.22}
\newcommand\theotherbound{1.182}
\newcommand\isomplus{\mathop{{\rm Isom}_+}}

\begin{document}

\begin{asciiabstract}
 If a closed, orientable hyperbolic 3--manifold M has volume at most
 1.22 then H_1(M;Z_p) has dimension at most 2 for every prime p not 2
 or 7, and H_1(M;Z_2) and H_1(M;Z_7) have dimension at most 3.  The
 proof combines several deep results about hyperbolic 3--manifolds.
 The strategy is to compare the volume of a tube about a shortest
 closed geodesic C in M with the volumes of tubes about short closed
 geodesics in a sequence of hyperbolic manifolds obtained from M by
 Dehn surgeries on C.
\end{asciiabstract}

\begin{htmlabstract}
If a closed, orientable hyperbolic 3&ndash;manifold M has volume at
most 1.22 then H<sub>1</sub>(M; <b>Z<sub></b>p</sub>) has dimension at
most 2 for every prime p not 2 or 7, and H<sub>1</sub>(M; 
<b>Z</b><sub>2</sub>) and H<sub>1</sub>(M; <b>Z</b><sub>7</sub>) have
dimension at most 3.  The proof combines several deep results about
hyperbolic 3&ndash;manifolds.  The strategy is to compare the volume
of a tube about a shortest closed geodesic C&sub; M with the volumes
of tubes about short closed geodesics in a sequence of hyperbolic
manifolds obtained from M by Dehn surgeries on C.
\end{htmlabstract}

\begin{abstract} 
  If a closed, orientable hyperbolic $3$--manifold $M$ has volume at
  most 1.22 then $H_1(M;\mathbb Z_p)$ has dimension at most $2$ for
  every prime $p\ne2,7$, and $H_1(M;\mathbb Z_2)$ and $H_1(M;\mathbb Z_7)$
  have dimension at most $3$.  The proof combines several deep results
  about hyperbolic $3$--manifolds.  The strategy is to compare the
  volume of a tube about a shortest closed geodesic $C\subset M$ with
  the volumes of tubes about short closed geodesics in a sequence of
  hyperbolic manifolds obtained from $M$ by Dehn surgeries on $C$.

\end{abstract}

\maketitle

\section{Introduction}

We shall prove:

\Theorem\label{main theorem} Suppose that $M$ is a closed, orientable
hyperbolic $3$--manifold with volume at most $\thebound$. Then
$H_1(M;\ZZ_p)$ has dimension at most $2$ for every prime $p\ne2,7$,
and $H_1(M;\ZZ_2)$ and $H_1(M;\ZZ_7)$ have dimension at most $3$.
Furthermore, if $M$ has volume at most $\theotherbound$, then
$H_1(M;\ZZ_7)$ has dimension at most $2$.  \EndTheorem

The bound of $2$ for the dimension of $H_1(M;\ZZ_p)$ is sharp when $p$
is $3$ or $5$.  Indeed, the manifolds {\tt m003(-3,1)}, and {\tt
  m007(3,1)} from the list given in \cite{snappea} have respective
volumes $ 0.94\ldots$ and $1.01\ldots$, and their integer homology
groups are respectively isomorphic to $\ZZ_5 \oplus \ZZ_5$ and
$\ZZ_3\oplus\ZZ_6$.

Apart from these two examples, the only example known to us of a
closed, orientable hyperbolic $3$--manifold with volume at most
$\thebound$ is the manifold {\tt m003(-2,3)} from the list given in
\cite{snappea}. These three examples suggest that the bounds for the
dimension of $H_1(M;\ZZ_p)$ given by \fullref{main theorem} may
not be sharp for $p\ne3,5$.

The proof of \fullref{main theorem} depends on several deep results,
including a strong form of the ``$\log3$ Theorem'' of Anderson,
Canary, Culler and Shalen \cite{ACCS,CS}; the Embedded Tube Theorem of
Gabai, Meyerhoff and N Thurston \cite{GMT}; the Marden Tameness
Conjecture, recently proved by Agol \cite{Agol} and by Calegari and
Gabai \cite{Calegari-Gabai}; and an even more recent result due to
Agol, Dunfield, Storm and W Thurston \cite{ADST}. The strategy of our
proof is to compare the volume of a tube about a shortest closed
geodesic $C\subset M$ with the volumes of tubes about short closed
geodesics in a sequence of hyperbolic manifolds obtained from $M$ by
Dehn surgeries on $C$.

After establishing some basic conventions in \fullref{convention
  section}, we carry out the strategy described above in Sections
\ref{drilling section}--\ref{odd section}, for the case of manifolds
which are ``non-exceptional'' in the sense that they contain shortest
geodesics with tube radius greater than $(\log3)/2$. In \fullref{arbitrary section}, for the case of non-exceptional manifolds
with volume at most $\thebound$, we establish a bound of $3$ for the
dimension of $H_1(M;\ZZ_p)$ for any prime $p$.  In \fullref{odd
  section}, again for the case of non-exceptional manifolds with
volume at most $\thebound$, we establish a bound of $2$ for the
dimension of $H_1(M;\ZZ_p)$ for any odd prime $p$. In \fullref{exceptional section} we use results from \cite{GMT} to handle the
case of exceptional manifolds, and complete the proof of \fullref{main theorem}.

The research described in this paper was partially supported by NSF
grants DMS-0204142 and DMS-0504975.

\section{Definitions and conventions} \label{convention section}
\Number If $g$ is a loxodromic isometry of hyperbolic $3$--space
$\HH^3$ we shall let $A_g$ denote the hyperbolic geodesic which is the
axis of $g$.  The {\it cylinder about $A_g$ of radius $r$} is the open
set $Z_r(g) = \{ x\in \HH^3 \;\vert\; \dist(x,A_g) < r \}$.
\EndNumber

\Number\label{number one} Suppose that $M$ is a complete, {orientable}
hyperbolic $3$--manifold.  Let us identify $M$ with $\HH^3/\Gamma$,
where $\Gamma\cong\pi_1(M)$ is a discrete, torsion-free subgroup of
$\isomplus\HH^3$.  If $C$ is a simple closed geodesic in $M$ then
there is a loxodromic isometry $g\in\Gamma$ with $A_g/\langle g\rangle
= C$.  For any $r>0$ the image $Z_{r}(g)/\langle g \rangle$ of
$Z_r(g)$ under the covering projection is a neighborhood of $C$ in
$M$.  For sufficiently small $r>0$ we have
  $$\{h\in\Gamma \;|\; h(Z_r(g))\cap Z_r(g)\not=\emptyset\} = \langle
  g\rangle .$$
Let $R$ denote the supremum of the set of $r$ for which this condition
holds.  We define $\tube(C) = Z_{R}(g)/\langle g \rangle$ to be the
{\it maximal tube about $C$}. We shall refer to $R$ as the {\it
tube radius} of $C$, and denote it by $\tuberad(C)$.
\EndNumber

\Number \label{number two} If $C$ is a simple closed geodesic in a
closed hyperbolic $3$--manifold $M$, it follows from \cite{kojima},
\cite{drilling} that $M-C$ is homeomorphic to a hyperbolic manifold
$N$ of finite volume having one cusp.  The manifold $N$, which by
Mostow rigidity is unique up to isometry, will be denoted
$\drillc(M)$.  \EndNumber

\Number If $C$ is a shortest closed geodesic in a closed hyperbolic
$3$--manifold $M$, ie, one such that $\length (C)\le\length
(C')$ for every other closed geodesic $C'$, then in particular $C$ is
simple, and the notions of \ref{number one} and \ref{number two} apply
to $C$.  \EndNumber

\Number Suppose that $N=\HH^3/\Gamma$ is a non-compact orientable
complete hyperbolic manifold of finite volume.  Let
$\Pi\cong\ZZ\times\ZZ$ be a maximal parabolic subgroup of $\Gamma$ (so
that $\Pi$ corresponds to a peripheral subgroup under the isomorphism
of $\Gamma$ with $\pi_1(N)$).  Let $\xi$ denote the fixed point of
$\Pi$ on the sphere at infinity and let $B$ be an open horoball
centered at $\xi$ such that $\{g\in \Gamma \;|\;gB\cap B
\not=\emptyset\} = \Pi$.  Then $\calh = B/\Pi$, which we identify with
the image of $B$ in $N$, is called a {\it cusp neighborhood} in $N$.

If $\calh$ is a cusp neighborhood in $N=\HH^3/\Gamma$ then the inverse
image of $\calh$ under the covering projection $\HH^3\to N$ is a union
of disjoint open horoballs.  The cusp neighborhood $\calh$ is maximal
if and only there exist two of these disjoint horoballs whose closures
have non-empty intersection.  \EndNumber

\Number If $N$ is a complete, orientable hyperbolic manifold of finite
volume, $\hat N$ will denote a compact core of $N$. Thus $\hat N$ is a
compact $3$--manifold whose boundary components are all tori, and the
number of these tori is equal to the number of cusps~of~$N$.
\EndNumber

\section{Drilling and packing}\label{drilling section}

\Lemma \label{agoldunfield} Suppose that $M$ is a closed, orientable
hyperbolic $3$--manifold, and that $C$ is a shortest geodesic in $M$.
Set $N=\drillc(M)$. If $\tuberad(C) \ge (\log3)/2$ then $\vol N <
3.0177\vol M$.  \EndLemma

\Proof The proof is based on a result due to Agol, Dunfield, Storm and W
Thurston \cite{ADST}.  We let $L$ denote the length of the geodesic
$C$ in the closed hyperbolic $3$--manifold $M$, and we set $R=\tuberad
(C)$ and $T=\tube(C)$.  Proposition 10.1 of \cite{ADST} states that
\begin{align*}
\vol N &\le (\coth^3 2R)(\vol M + \frac{\pi}{2}L\tanh R \tanh 2R) .
\\
\tag*{\hbox{Note that}}
\vol T &= \pi L\sinh^2R 
= \left(\frac{\pi}{2}L\tanh R\right)(2\sinh R \cosh R)
\\
&\phantom{=\pi L\sinh^2R }\,=  \left(\frac{\pi}{2}L\tanh R\right)(\sinh 2R) .
\\
\tag*{\hbox{Thus}}
\vol N &\le (\coth^3 2R)\left(\vol M + \vol T\frac{\tanh 2R}{\sinh 2R}\right)\\
       &= (\coth^3 2R)\left(\vol M + \frac{\vol T}{\cosh 2R}\right).\\
\end{align*}
In the language of \cite{prez}, the quantity $(\vol T)/(\vol M)$ is
the density of a tube packing in $\HH^3$. According to \cite[Corollary
4.4]{prez}, we have $(\vol T)/\vol M < 0.91$.  Hence $\vol N <
f(x)\vol(M)$, where $f(x)$ is defined for $x\ge0$ by
$$
f(x) = (\coth^3 2x)\left(1 + \frac{0.91}{\cosh 2x}\right).$$
Since
$f(x)$ is decreasing for $x\ge 0$, and since a direct computation
shows that $f(0.5495)=3.01762\ldots$, we have $\vol N<3.0177\vol M$
whenever $R\ge0.5495$.

It remains to consider the case in which
$0.5495>R\ge(\log3)/2=0.5493\ldots$. In this case we use \cite[Theorem
4.3]{prez}, which asserts that the tube-packing density $(\vol T)/\vol
M$ is bounded above by $(\sinh R)g(R)$, where $g(x)$ is defined for
$x>0$ by
$$g(x)=\frac{\arcsin\frac{1}{2\cosh r}}{\arcsinh\frac{\tanh
    r}{\sqrt3}}.$$
Since $g(x)$ is clearly a decreasing function for $x>0$, and since
$\sinh R$ is increasing for $x>0$, we have
$$(\vol T)/(\vol M)<(\sinh0.5495)g((\log3)/2)=0.90817\ldots.$$
Hence 
$\vol N < f_1(x)\vol(M)$, where $f_1(x)$ is defined for $x\ge0$ by
$$
f_1(x) = (\coth^3 2x)\left(1 + \frac{0.90817}{\cosh 2x}\right).$$
Again, $f_1(x)$ is decreasing for $x\ge 0$, and we see by direct
computation that\break $f_1((\log3)/2)=3.017392\ldots$. Hence we have $\vol
N<3.0174\vol M$ in this case.  \EndProof

\Lemma\label{volume bound for large R} Suppose that $M$ is a closed,
orientable hyperbolic $3$--manifold such that $\vol M \le \thebound$,
and that $C$ is a shortest geodesic in $M$.  Set $N=\drillc(M)$.  If
$\tuberad(C) > (\log3)/2$ then the maximal cusp neighborhood in $N$
has volume less than $\pi$.  \EndLemma

\Proof We let $d(\infty) = .853276\ldots$ denote B\"or\"oczky's lower
bound \cite{bor} for the density of a horoball packing in hyperbolic
space.  It follows from the definition of the density of a horoball
packing that the volume of a maximal cusp neighborhood in $N$ is at
most $d(\infty)\vol N$.  \fullref{agoldunfield} gives $\vol N <
3.0177\cdot\thebound<\pi/d(\infty)$, and the conclusion follows.
\EndProof

\section{Filling}

As in \cite{ACCS}, we shall say that a group is {\it semifree} if it
is a free product of free abelian groups; and we shall say that a
group $\Gamma$ is {\it $k$--semifree} if every subgroup of $\Gamma$
whose rank is at most $k$ is semifree. Note that $\Gamma$ is
$2$--semifree if and only if every rank-$2$ subgroup of $\Gamma$ is
either free or free abelian.

The following improved version of \cite[Theorem 6.1]{ACCS} is made
possible by more recent developments.

\Theorem\label{better 6.1} Let $k\ge 2$ be an integer and let $ \F $
be a Kleinian group which is freely generated by elements
$\xi_1,\ldots,\xi_k$.  Let $z$ be any point of $ \HH^3$ and set
$d_i=\dist(z,\xi_i\cdot z)$ for $i=1,\ldots,k$. Then we have
$$\sum_{i=1}^k \frac{1}{1+e^{d_i}}\leq\frac12.$$
In particular there
is some $i\in\{1,\ldots,k\}$ such that $d_i\geq\log(2k-1)$.
\end{theorem}

\Proof If $\Gamma$ is geometrically finite this is included in
\cite[Theorem 6.1]{ACCS}. In the general case, $\Gamma$ is
topologically tame according to \cite{Agol} and \cite{Calegari-Gabai},
and it then follows from \cite[Theorem 1.1]{Oh}, or from the
corresponding result for the free case in \cite{NS}, that $\Gamma$ is
an algebraic limit of geometrically finite groups; more precisely,
there is a sequence of geometrically finite Kleinian groups
$(\Gamma_j)_{j\ge1}$ such that each $\Gamma_j$ is freely generated by
elements $\xi_{1j},\ldots,\xi_{kj}$, and $\lim_{j\to\infty}
\xi_{ij}=\xi_i$ for $i=1,\ldots,k$. Given any $z\in\HH^3$, we set
$d_{ij}=\dist(z,\xi_{ij}\cdot z)$ for each $j\ge1$ and for
$i=1,\ldots,k$. According to \cite[Theorem 6.1]{ACCS}, we have
$$\sum_{i=1}^k
\frac{1}{1+e^{d_{ij}}}\leq\frac12$$
for each $j\ge1$. Taking limits as $j\to\infty$ we conclude that 
$$\sum_{i=1}^k
\frac{1}{1+e^{d_{i}}}\leq\frac12.\proved$$
\EndProof

Let us also recall the following definition from \cite[Section
8]{ACCS}.  Let $\Gamma$ be a discrete torsion-free subgroup of
$\Isom_+(\HH^3)$.  A positive number $\lambda$ is termed a {\em strong
  Margulis number} for $\Gamma$, or for the orientable hyperbolic
3--manifold $N=\HH^3/\Gamma$, if whenever $\xi $ and $\eta$ are
non-commuting elements of $\Gamma$, we have $$\frac{1}{1+e^{\dist(\xi
    \cdot z,z)}} + \frac{1}{1+e^{\dist(\eta\cdot z,z)}}\le
\frac{2}{1+e^{\lambda}}\ .$$

The following improved version of \cite[Proposition 8.4]{ACCS} is an
immediate consequence of \fullref{better 6.1}.

\Corollary\label{better 8.4} Let $\Gamma$ be a discrete subgroup of
$\Isom_+(\HH^3)$. Suppose that $\Gamma$ is $2$--semifree.  Then $\log3$
is a strong Margulis number for $\Gamma$. \NoProof \EndCorollary

\Lemma\label{big fillings} Let $N$ be a non-compact finite-volume
hyperbolic $3$--manifold.  Suppose that $S$ is a boundary component of
the compact core $\hat N$, and $\calh$ is the maximal cusp
neighborhood in $N$ corresponding to $S$.  If infinitely many of the
manifolds obtained by Dehn filling $\hat N$ along $S$ have
$2$--semifree fundamental group then $\calh$ has volume at least $\pi$.
\EndLemma

\Proof Suppose that $(N_i)$ is an infinite sequence of distinct
hyperbolic manifolds obtained by Dehn filling $\hat N$ along $S$, and
that $\pi_1(N_i)$ is $2$--semifree for each $i$.

Thurston's Dehn filling theorem \cite[Appendix B]{BP}, implies that
for each sufficiently large $i$, the manifold $N_i$ admits a
hyperbolic metric; that the core curve of the Dehn filling $N_i$ of
$\hat N$ is isotopic to a geodesic $C_i$ in $N_i$; that the length
$L_i$ of $C_i$ tends to $0$ as $i\to\infty$; and that the sequence of
maximal tubes $(\tube(C_i))_{i\ge1}$ converges geometrically to
$\calh$.  In particular
$$\lim_{i\to\infty}\vol (\tube(C_i)) = \vol\calh .$$

According to \fullref{better 8.4}, $\log 3$ is a strong Margulis
number for each of the hyperbolic manifolds $N_i$.  It therefore
follows from \cite[ Corollary 10.5]{ACCS} that $\vol \tube(C_i) >
V(L_i)$, where $V$ is an explicitly defined function such that
$\lim_{x\to0}V(x)=\pi$.  In particular, this shows that
$$\vol \calh \ge \lim_{i\to\infty}V(L_i)\ge\pi.\proved$$
\EndProof

\section{Non-exceptional manifolds, arbitrary primes}\label{arbitrary section}

\Number A closed hyperbolic $3$--manifold $M$ will be termed {\it
  exceptional} if every shortest geodesic in $M$ has tube radius at
most $(\log 3)/2$.  \EndNumber

In this section we shall prove a result, \fullref{large R},
which gives a bound of $3$ for the dimension of $H_1(M;\ZZ_p)$ for any
prime $p$ when $M$ is a non-exceptional manifold with volume at most
$\thebound$.
\eject
\Lemma\label{no finite index} Suppose that $M$ is a compact,
irreducible, orientable $3$--manifold, such that every non-cyclic
abelian subgroup of $\pi_1(M)$ is carried by a torus component of
$\partial M$. Suppose that either\begin{itemize}
\item[\rm(i)] $\dim H_1(M;\QQ)\ge3$, or \item[\rm(ii)]
$M$ is closed and $\dim H_1(M;\ZZ_p)\ge4$ for some prime
$p$. \end{itemize}Then
$\pi_1(M)$ is $2$--semifree.  \EndLemma

\Proof Let $X$ be any subgroup of $\pi_1(M)$ having rank at most $2$.
According to \cite[Theorem VI.4.1]{JS}, $X$ is free, or free abelian,
or of finite index in $\pi_1(M)$. If $\dim H_1(M;\QQ)\ge3$, it is
clear that $X$ has infinite index in $\pi_1(M)$. If $M$ is closed and
$H_1(M;\ZZ_p)\ge4$ for some prime $p$, then Proposition 1.1 of
\cite{shalen-wagreich} implies that every $2$--generator subgroup of
$\pi_1(M)$ has infinite index. Thus in either case $X$ is either free
or free abelian. This shows that $\pi_1(M)$ is $2$--semifree.
\EndProof

\Proposition\label{large R} Suppose that $M$ is a closed, orientable,
non-exceptional hyperbolic $3$--manifold such that $\vol M \le
\thebound$.  Then $H_1(M;\ZZ_p)$ has dimension at most $3$ for every
prime $p$.  \EndProposition

\Proof Since $M$ is non-exceptional, there is a shortest geodesic $C$
in $M$ with $R=\tuberad(C) > (\log3)/2$. We set $N=\drillc (M)$. Let
$\calh$ denote the maximal cusp neighborhood in $N$.  Since $R >
(\log3)/2$, \fullref{volume bound for large R} implies that
$\vol \calh < \pi$.

Now assume that $\dim H_1(M;\ZZ_p)\ge4$ for some prime $p$.  There is
an infinite sequence $(M_i)$ of manifolds obtained by distinct Dehn
fillings of $\hat N$ such that $H_1(M_i;\ZZ_p)$ has dimension at least
$4$ for each $i$.  (For example, if $(\lambda,\mu)$ is a basis for
$H_1(\partial \hat N,\ZZ_p)$ such that $\lambda$ belongs to the kernel
of the inclusion homomorphism $H_1(\partial \hat N,\ZZ_p)\to H_1(\hat
N,\ZZ_p)$, we may take $M_i$ to be obtained by the Dehn surgery
corresponding to a simple closed curve in $\partial \hat N$
representing the homology class $\lambda+ip\mu$.)  It follows from
Thurston's Dehn filling theorem \cite[Appendix B]{BP} that for
sufficiently large $i$ the manifold $M_i$ is hyperbolic. Hence by case
(ii) of \fullref{no finite index}, the fundamental group of $M_i$ is
$2$--semifree for sufficiently large $i$.  Thus \fullref{big
  fillings} implies that $\vol \calh\ge\pi$, a contradiction.
\EndProof

\section{Non-exceptional manifolds, odd primes}\label{odd section}

\fullref{other main theorem}, which is proved in this section,
gives a bound of $2$ for the dimension of $H_1(M;\ZZ_p)$ for any odd
prime $p$ when $M$ is a non-exceptional manifold with volume at most
$\thebound$.

\Definition Let $N$ be a connected manifold, $\star\in N$ a base
point, and $Q$ a subgroup of $\pi_1(N,\star)$. We shall say that a
connected based covering space $r:(N',\star')\to (N,\star)$ {\it
  carries} the subgroup $Q$ if $Q\le r_\sharp(\pi_1(
N',\star'))\le\pi_1(N,\star)$ \EndDefinition

\Lemma\label{beta} Suppose that $\calh$ is a maximal cusp neighborhood
in a finite-volume hyperbolic $3$--manifold $N$. Let $\star$ be a base
point in $\calh$, and let $P\le\pi_1(N,\star)$ denote the image of
$\pi_1(\calh,\star)$ under inclusion. Then there is an element $\beta$
of $\pi_1(N,\star)$ with the following property: \Parts
\Part{($\dagger$)}For every based covering space $r:(
N',\star')\to (N,\star)$ which carries the subgroup $\langle
P,\beta\rangle$ of $\pi_1(N,\star)$, there is a maximal cusp
neighborhood $\calh'$ in $N'$ which is isometric to $\calh$.
\EndParts
\EndLemma

\Proof. We write $N=\HH^3/\Gamma$, where $\Gamma$ is a discrete,
torsion-free subgroup of $\Isom(\HH^3)$. Let $q:\HH^3\to N$ denote the
quotient map and fix a base point $\star'$ which is mapped to $\star$
by $q$. The components of $q^{-1}(\calh)$ are horoballs.  Let $B_0$
denote the component of $q^{-1}(\calh)$ containing $\star'$.  The
stabilizer $\Gamma_0$ of $B_0$ is mapped onto the subgroup $P$ of
$\pi_1(N,\star)$ by the natural isomorphism
$\iota:\Gamma\to\pi_1(N,\star)$.

Since $\calh$ is a maximal cusp, there is a component $B_1\ne B_0$ of
$q^{-1}(\calh)$ such that $\overline
{B_1}\cap\overline{B_0}\not=\emptyset$.  We fix an element $g$ of
$\Gamma$ such that $g(B_0)=B_1$, and we set $\beta=\iota(g)\in
\pi_1(N,\star)$.

To show that $\beta$ has property ($\dagger$), we consider an
arbitrary based covering space $r:(N',\star')\to (N,\star)$ which
carries the subgroup $\langle P,\beta\rangle$ of $\pi_1(N,\star)$. We
may identify $N'$ with $\HH^3/\Gamma'$, where $\Gamma'$ is some
subgroup of $\Gamma$ containing $\langle\Gamma_0,g\rangle$.

Since $\Gamma_0\subset\Gamma'$, the cusp neighborhood $\calh$ lifts to
a cusp neighborhood $\calh'$ in $N'$. In particular $\calh'$ is
isometric to $\calh$. The horoballs $B_0$ and $B_1=g(B_0)$ are
distinct components of $ (q')^{-1}(\calh')$, where $ q':\HH^3\to N'$
denotes the quotient map. Since $g\in\Gamma'$ and $\overline
{B_1}\cap\overline{B_0}\not=\emptyset$, the cusp neighborhood $
\calh'$ is maximal.  \EndProof

\Proposition\label{other main theorem} Suppose that $M$ is a closed,
orientable, non-exceptional hyperbolic $3$--manifold such that $\vol M
\le \thebound$. Then $H_1(M;\ZZ_p)$ has dimension at most $2$ for
every odd prime $p$.  \EndProposition

\Proof Since $M$ is non-exceptional, we may fix a shortest geodesic
$C$ in $M$ with $R=\tuberad(C) > (\log3)/2$. We set $N=\drillc (M)$.
Let $\calh$ denote the maximal cusp neighborhood in $N$.  Since $R >
(\log3)/2$, \fullref{volume bound for large R} implies that
$\vol \calh < \pi$.

As in the statement of \fullref{beta}, we fix a base point $\star\in
\calh$, and we denote by $P\le\pi_1(N,\star)$ the image of
$\pi_1(\calh,\star)$ under inclusion. We fix an element $\beta$ of
$\pi_1(N,\star)$ having property ($\dagger$) of \fullref{beta}. We
set $Q=\langle P,\beta\rangle\le\pi_1(N,\star)$.

Suppose that $\dim H_1(M;\ZZ_p)\ge3$ for some prime $p$. We shall
prove the proposition by showing that this assumption leads to a
contradiction if $p$ is odd.

It follows from Poincar\'e duality that the image of the inclusion
homomorphism $\alpha:H_1(\partial \hat N;\ZZ_p)\to H_1(\hat N;\ZZ_p)$
has rank $1$.  Hence the image of $P$ under the natural homomorphism
$\pi_1(N,\star)\to H_1(N;\ZZ_p)$ has dimension $1$. It follows that
the image $\bar Q$ of $Q$ under this homomorphism has dimension either
$1$ or $2$. In the case $\dim\bar Q=1$ we shall obtain a contradiction
for any prime $p$.  In the case $\dim\bar Q=2$ we shall obtain a
contradiction for any odd prime $p$.

First consider the case $\dim\bar Q=1$.  We have assumed 
$\dim H_1(M;\ZZ_p)\ge3$.  Thus there is a $\ZZ_p\times\ZZ_p$--regular based
covering space $(N',\star')$ of $(N,\star)$ which carries $Q$. By
property ($\dagger$), there is a maximal cusp neighborhood $\calh'$ in
$N'$ which is isometric to $\calh$. In particular $\vol \calh' < \pi$.

Since in particular $(N',\star')$ carries $P$, the boundary of the
compact core $\hat N$ lifts to $\hat N'$. As $N'$ is a $p^2$--fold
regular covering, it follows that $\hat N'$ has $p^2\ge4$ boundary
components.  

It follows from Thurston's Dehn filling theorem \cite[Appendix B]{BP}
that there are infinitely many hyperbolic manifolds obtained by Dehn
filling one boundary component of $\hat N'$.  If $Z$ is any hyperbolic
manifold obtained by such a filling, then $Z$ has at least three
boundary components, and it follows from case (i) of \fullref{no
  finite index} that $\pi_1(Z)$ is $2$--semifree. It therefore follows
from \fullref{big fillings} that each maximal cusp neighborhood in
$N'$ has volume at least $\pi$. Since we have seen that $\vol \calh' <
\pi$, this gives the desired contradiction in the case $\dim\bar Q=1$.

It remains to consider the case in which $\dim\bar Q=2$ and the prime
$p$ is odd. Since we have assumed that $\dim H_1(M;\ZZ_p)\ge3$, there
is a $p$--fold cyclic based covering space $(N',\star')$ of $(N,\star)$
which carries $Q$.  Since $N'$ carries $P$, the boundary of the
compact core $\hat N$ lifts to $\hat N'$, and as $N'$ is a $p$--fold
regular covering, it follows that $\hat N'$ has $p$ boundary
components.
  
We claim that the inclusion homomorphism $\alpha':H_1(\partial\hat
N',\ZZ_p)\to H_1(\hat N',\ZZ_p)$ is not surjective. To establish this,
we consider the commutative diagram
$$
\xymatrix{
H_1(\partial \hat N';\ZZ_p) \ar@{->}^{\alpha'}[r] \ar@{->}[d] &
   H_1(N';\ZZ_p) \ar@{->}[d]^{r_*}\\
H_1(\partial \hat N;\ZZ_p) \ar@{->}_\alpha[r]& H_1(N;\ZZ_p)\\
}
$$
where $r:N'\to N$ is the covering projection. Since $(N',\star')$
carries $Q$ we have $\bar Q\subset \image r_*$. Hence surjectivity of
$\alpha'$ would imply $\bar Q\subset \image\alpha$. This is
impossible: we observed above that $\image\alpha$ has rank $1$, and we
are in the case $\dim\bar Q=2$. Thus $\alpha'$ cannot be surjective.

Since $\hat N'$ has $p$ boundary components, it follows from
Poincar\'e duality that\break $\dim\image\alpha'=p\ge3$.  Since $\alpha'$ is
not surjective and $p$ is an odd prime, it follows that $\dim
H_1(N';\ZZ_p)\ge p+1\ge4$.

Since $(N',\star')$ carries $Q$, some subgroup $Q'$ of
$\pi_1(N',\star')$ is mapped isomorphically to $Q$ by $r_\sharp$.  In
particular $Q'$ has rank at most $3$.  Since $\dim H_1(N';\ZZ_p)\ge
4$, there is a $p$--fold cyclic based covering space $(N'',\star'')$ of
$(N',\star')$ which carries $Q'$. Hence $(N'',\star'')$ is a
$p^2$--fold (possibly irregular) based covering space of $(N,\star)$
which carries $Q$.  By property ($\dagger$), there is a maximal cusp
neighborhood $\calh''$ in $N''$ which is isometric to $\calh$. In
particular $\vol \calh'' < \pi$.

Since $P\le Q$, there is a component $T$ of $\partial\hat N'$ such
that $Q'$ contains a conjugate of the image of $\pi_1(T)$ under the
inclusion homomorphism $\pi_1(T)\to\pi_1(N')$. Hence $T$ lifts to the
$p$--fold cyclic covering space $N''$ of $N'$. It follows that the
covering projection $r':N''\to N'$ maps $p\ge3$ components of
$(r')^{-1}(\partial \hat N')$ to $T$. As $\hat N'$ has at least three
boundary components, $\hat N''$ must have at least five boundary
components.

Hence if $Z$ is any hyperbolic manifold obtained by Dehn filling one
boundary component of $\hat N''$, we have $\dim H_1(Z;\QQ)\ge 4>3$,
and it follows from case (i) of \fullref{no finite index} that
$\pi_1(Z)$ is $2$--semifree. It therefore follows from \fullref{big
  fillings} and Thurston's Dehn filling theorem that each maximal cusp
neighborhood in $N''$ has volume at least $\pi$. Since we have seen
that $\vol \calh'' < \pi$, we have the desired contradiction in this
case as well.  \EndProof

\section{Exceptional manifolds}\label{exceptional section}

Our treatment of exceptional manifolds begins with \fullref{GMT stuff} below, the proof of which will largely consist of
citing material from \cite{GMT}. In order to state it we must first
introduce some notation.

For $k=0,\ldots, 6$ we define constants $\tau_k$ as follows:

{\parindent=1in\parskip=0pt
$\tau_0 =  0.4779$

$\tau_1 =  1.0756$

$\tau_2 =  1.0527$

$\tau_3 =  1.2599$

$\tau_4 =  1.2521$

$\tau_5 =  1.0239$

$\tau_6 =  1.0239$
}

For $k=0,\ldots,6$ let $\ek$ be the $2$--generator group with
presentation 
$$\ek = \langle x,y \colon r_{1,k}, r_{2,k}\rangle ,$$
where the relators
$r_{1,k}=r_{1,k}(x,y)$ and $r_{2,k}=r_{2,k}(x,y)$ are the words listed
below (in which we have set $X=x^{-1}$ and $Y=y^{-1}$):

{\parindent=1in\parskip=0pt
$r_{1,0} = xyXyyXyxyy$,

$r_{2,0} = XyxyxYxyxy$,

\medskip\goodbreak
$r_{1,1} = XXyXYXYxYXYXyXXyy$,

$r_{2,1} = XXyyXyxyxYxyxyXyy$,

\medskip\goodbreak
$r_{1,2} = XyxyxYxxYxyxyXyy$,

$r_{2,2} = XXyXXyyXyxyxyXyy$,

\medskip\goodbreak
$r_{1,3} = XXyxyXXyyXYXyXYxYXYxxYXYxYXyXYXyy$,

$r_{2,3} = XXyxyXyxYxyxYYxyxYxyXyxyXXyyXYXyy$,

\medskip\goodbreak
$r_{1,4} = XXyxyXyxYxyxYxyXyxyXXyyXYXyXYXyy$,

$r_{2,4} = XXyxyXyxyXXyyXYXyXYxYXYxYXyXYXyy$,

\medskip\goodbreak
$r_{1,5} = XyXYXyXyxyxYxyxy$,

$r_{2,5} = XyxyxYxYXYxYxyxy$,

\medskip\goodbreak
$r_{1,6} = XYXyXYxYXyXYXyxy$,

$r_{2,6} = XYXyxyXyxYxyXyxy$.
}

The group $\ee_0$ is the fundamental group of an arithmetic hyperbolic
$3$--manifold which is known as $\volthree$.  This manifold, which was
studied in \cite{jones-reid}, is described as {\tt m007(3,1)} in the
list given in \cite{snappea}, and can also be described as the
manifold obtained by a $(-1,2)$ Dehn filling of the once-punctured
torus bundle with monodromy $-R^2L$.

\Proposition\label{GMT stuff} Suppose that $M$ is an exceptional
closed, orientable hyperbolic $3$--manifold which is not isometric to
$\volthree$. Then there exists an integer $k$ with $1\le k\le6$ such
that the following conditions hold: \Parts
\Part{(1)}$M$ has a finite-sheeted cover $\widetilde M$ such that
$\pi_1(\widetilde M)$ is isomorphic to a quotient of $\ek$; and
\Part{(2)}there is a shortest closed geodesic $C$ in $M$ such that
$\vol(\tube(C))\ge\tau_k$.
\EndParts
\EndProposition

\Proof This is in large part an application of results from
\cite{GMT}, and we begin by reviewing some material from that paper.

We begin by considering an arbitrary simple closed geodesic $C$ in a
closed, orientable hyperbolic $3$--manifold $M=\HH^3/\Gamma$. As we
pointed out in \ref{number one}, there is a loxodromic isometry
$\eff\in\Gamma$ with $A_\eff/\langle\eff\rangle=C$. If we set
$R=\tuberad(C)$ and $Z=Z_R(\eff)$, it follows from the definitions
that $\tube(C)=Z/\langle \eff\rangle$, that $h(Z)\cap Z=\emptyset$ for
every $h\in\Gamma-\langle \eff\rangle$, and that there is an element
$\omigosh\in\Gamma-\langle \eff\rangle$ such that $\omigosh(\bar
Z)\cap \bar Z\ne\emptyset$.

Let us define an ordered pair $(\eff,\omigosh)$ of elements of
$\Gamma$ to be a {\it \gimpty pair} for the simple geodesic $C$ if we
have (i) $A_\eff/\langle\eff\rangle=C$, (ii) $\omigosh\notin\langle
\eff\rangle$, and (iii) $\omigosh(\bar Z)\cap \bar Z\ne\emptyset$.
Note that since $\langle\eff\rangle$ must be a maximal cyclic subgroup
of $\Gamma$, condition (ii) implies that the group
$\langle\eff,\omigosh\rangle$ is non-elementary.
 
Set $\calp=\{(L,D,R)\in \CC^3:\Re L,\Re D>0\}$.  For any point
$\point=(L,D,R)\in\calp$ we will denote by
$(\eff_\point,\omigosh_\point)$ the pair
$(\eff,\omigosh)\in\Isom_+(\HH^3)\times\Isom_+(\HH^3)$, where
$\eff,\omigosh\in PGL_2(\CC)=\Isom_+(\HH^3)$ are defined by
\begin{gather*}
\eff=\bigg[\begin{matrix}e^{L/2}&0\cr0&e^{-L/2}\end{matrix}\bigg]\\
\tag*{\hbox{and}}
\omigosh=\bigg[\begin{matrix}e^{R/2}&0\cr0&e^{-R/2}\end{matrix}\bigg]\ 
\bigg[\begin{matrix}1&\hfill1\cr1&-1\end{matrix}\bigg]\ 
\bigg[\begin{matrix}e^{D/2}&0\cr0&e^{-D/2}\end{matrix}\bigg]\ 
\bigg[\begin{matrix}1&\hfill1\cr1&-1\end{matrix}\bigg]\ .
\end{gather*}
With this definition, $\eff_\point$ has (real) translation length $\Re L$, and
the (minimum) distance between $A_\eff$ and $\omigosh(A_\eff)$ is $(\Re
D)/2$.

In \cite[Section 1]{GMT}, it is shown that if $(\eff,\omigosh)$ is a
\gimpty pair for a shortest geodesic $C$ in a closed, orientable
hyperbolic $3$--manifold and $\tuberad(C)\le(\log3)/2$, then
$(\eff,\omigosh)$ is conjugate by some element of $\Isom^+(\HH^3)$ to
a pair of the form $(\eff_\point,\omigosh_\point)$ where
$\point\in\calp$ is a point such that
$\exp(\point)\doteq(e^L,e^D,e^R)$ lies in the union $X_0\cup\cdots\cup
X_6$ of seven disjoint open subsets of $\CC^3$ that are explicitly
defined in \cite[Proposition 1.28]{GMT}.

For every $k$ with $0\le k\le6$ and every point $\point=(L,D,R)$ such
that $\exp(\point)\in X_k$, it follows from \cite[Definition 1.27 and
Proposition 1.28]{GMT} that \Parts
\Part{(I)}
the isometries $r_{1,k}(\eff_\point,\omigosh_P)$
 and $r_{2,k}(\eff_\point,\omigosh_\point)$ have
 translation length  less than  $\Re L$; 
\EndParts
and it follows from \cite[Table 1.1]{GMT} that
\Parts
\Part{(II)}
$\pi\Re(L)\sinh^2(\Re(D)/2)>\tau_k$.  \EndParts

According to \cite[Proposition 3.1]{GMT}, if $C$ is a shortest
geodesic in a closed, orientable hyperbolic $3$--manifold, and if some
\gimpty pair for $C$ has the form $(\eff_\point,\omigosh_\point)$ for
some $\point$ with {$\exp(\point)\in X_0$}, then $M$ is isometric to
$\volthree$.
 
Now suppose that $M$ is an exceptional closed, orientable hyperbolic
$3$--manifold. Let us choose a shortest closed geodesic $C$ in $M$. By
the definition of an exceptional manifold, $C$ has tube radius
$\le(\log3)/2$. Hence the facts recalled above imply that $C$ has a
\gimpty pair of the form $(\eff_\point,\omigosh_\point)$ for some
$\point$ such that $\exp(\point)\in X_k$ for some $k$ with $0\le
k\le6$; and furthermore, that if $M$ is not isometric to $\volthree$,
then $1\le k\le6$. We shall show that conclusions (1) and (2) hold
with this choice of $k$.
 
For $i=1,2$ it follows from property (I) above that the element
$r_{i,k}(\eff,\omega)$ has real translation length less than the real
translation length $\Re L$ of $\eff$. Since $C$ is a shortest geodesic
in $M$, it follows that the conjugacy class of $r_{i,k}(\eff,\omega)$
is not represented by a closed geodesic in $M$. As $M$ is closed it
follows that $r_{i,k}(\eff,\omega)$ is the identity for $i=1,2$.
Hence the subgroup of $\Gamma$ generated by $\eff$ and $\omega$ is
isomorphic to a quotient of $\ek$. Since we observed above that
$\langle\eff,\omega\rangle$ is non-elementary, there is a non-abelian
subgroup $Y$ of $\pi_1(M)$ which is isomorphic to a quotient of $\ek$.
In particular $Y$ has rank $2$, and it cannot be a free group of rank
$2$ since the relators $r_{1,k}$ and $r_{2,k}$ are non-trivial. Hence
by \cite[Theorem VI.4.1]{JS} we must have $|\pi_1(M):Y|<\infty$. This
proves (1).

Finally, we recall that
$$\vol \tube(C)=\pi(\length (C))\sinh^2(\tuberad (C))=\pi (\Re L)\sinh^2((\Re D)/2).$$
Hence (2) follows from (II).
\EndProof

We shall also need the following slight refinement of
\cite[Proposition 1.1]{shalen-wagreich}.

\Proposition\label{refined SW} Let $p$ be a prime and let $M$ be a
closed $3$--manifold. If $p$ is odd assume that $M$ is orientable. Let
$X$ be a finitely generated subgroup of $\pi_1(M)$, and set $n=\dim
H_1(X;\ZZ_p)$.  If $\dim H_1(M;\ZZ_p)\ge \max(3,n+2)$, then $X$ has
infinite index in $\pi_1(M)$. In fact, $X$ is contained in infinitely
many distinct finite-index subgroups of $\pi_1(M)$.  \EndProposition

\Proof In this proof, as in \cite[Section 1]{shalen-wagreich}, for any
group $G$ we shall denote by $G_1$ the subgroup of $G$ generated by
all commutators and $p$-th powers, where $p$ is the prime given in the
hypothesis. Since $\dim H_1(X;\ZZ_p)=n$ we may write $X=EX_1$ for some
rank-$n$ subgroup $E$ of $X$.

We first assume that $n\ge 1$.  Set $\Gamma=\pi_1(M)$. Let $\cals$
denote the set of all finite-index subgroups $\Delta$ of $\Gamma$ such
that $\Delta\ge X$ and $\dim H_1(\Delta;\ZZ_p)\ge n+2$. The hypothesis
gives $\Gamma\in\cals$, so that $\cals\ne\emptyset$. Hence it suffices
to show that every subgroup $\Delta\in\cals$ has a proper subgroup $D$
such that $D\in\cals$.

Any group $\Delta\in\cals$ may be identified with $\pi_1(\widetilde
M)$ for some finite-sheeted covering space $\widetilde M$ of $M$. In
particular, $\widetilde M$ is a closed $3$--manifold, and is orientable
if $p$ is odd. Since $\Delta\in\cals$ we have $X\le\Delta=\pi_1(\tilde
M)$ and $\dim H_1(\widetilde M;\ZZ_p)=\dim H_1(\Delta;\ZZ_p)\ge n+2$.
Now set $D=E\Delta_1\le\Delta$.  Applying \cite[ Lemma
1.5]{shalen-wagreich}, with $\tilde M$ in place of $M$, we deduce that
$D$ is a proper, finite-index subgroup of $\Delta$, and that $\dim
H_1(D;\ZZ_p)\ge2n+1\ge n+2$. On the other hand, since
$\Delta\in\cals$, we have $X\le\Delta$, and hence $X=EX_1\le
E\Delta_1=D$. It now follows that $D\in\cals$, and the proof is
complete in the case $n\ge 1$.

If $n = 0$ then, since $\dim H_1(M;\ZZ_p)\ge 3$, there exists a
finitely generated subgroup $X'\ge X$ such that $H_1(X';\ZZ_p)$ has
dimension $1$.  The case of the Lemma which we have already proved
shows that $X'$ has infinite index.  Thus $X$ has infinite index as
well.  \EndProof

\Corollary\label{refined SW corollary} Let $p$ be a prime and let $M$
be a closed, orientable $3$--manifold. Let $X$ be a finite-index
subgroup of $\pi_1(M)$, and set $n=\dim H_1(X;\ZZ_p)$. Then $\dim
H_1(M;\ZZ_p)\le \max(2,n+1)$.  \NoProof\EndCorollary

\Lemma\label{small R} Suppose that $M$ is an exceptional hyperbolic
$3$--manifold with volume at most $\thebound$. Then $H_1(M;\ZZ_p)$ has
dimension at most $2$ for every prime $p\ne2,7$, and $H_1(M;\ZZ_2)$
and $H_1(M;\ZZ_7)$ have dimension at most $3$.  Furthermore, if $M$
has volume at most $\theotherbound$, then $H_1(M;\ZZ_7)$ has dimension
at most $2$.  \EndLemma

\Proof If $M$ is isometric to $\volthree$ then $\pi_1(M)$ is generated
by two elements, and the conclusions follow. For the rest of the proof
we assume that $M$ is not isometric to $\volthree$, and we fix an
integer $k$ with $1\le k\le6$ such that conditions (1) and (2) of
\fullref{GMT stuff} hold.

By condition (2) of \fullref{GMT stuff}, we may fix a shortest
closed geodesic $C$ in $M$ such that $\vol(T)\ge\tau_k$, where
$T=\tube(C)$.  It follows from a result of Przeworski's
\cite[Corollary 4.4]{prez} on the density of cylinder packings that
$\vol T < 0.91\vol M$, and so $\vol M>\tau_k/0.91$.  If $k=3$ we have
$\tau_k/0.91>\thebound$, and we get a contradiction to the hypothesis.
Hence $k\in\{1,2,4,5,6\}$.

Furthermore, we have $\tau_1/0.91>\theotherbound$. Hence if $\vol
M\le\theotherbound$ then $k\in\{2,4,5,6\}$.

By condition (1) of \fullref{GMT stuff}, $\pi_1(M)$ has a
finite-index subgroup $X$ which is isomorphic to a quotient of $\ek$.
From the defining presentations of the groups $\ee_1$, $\ee_2$,
$\ee_4$ $\ee_5$ and $\ee_6$, we find that $H_1(\ee_1;\ZZ)$ is
isomorphic to $\ZZ_7\oplus\ZZ_7$, that $H_1(\ee_2;\ZZ)$ and
$H_1(\ee_4;\ZZ)$ are isomorphic to $\ZZ_4\oplus\ZZ_{12}$, while
$H_1(\ee_5;\ZZ)$ and $H_1(\ee_6;\ZZ)$ are isomorphic to
$\ZZ_4\oplus\ZZ_{4}$. (One can check that the two groups $\ee_5$ and
$\ee_6$ are isomorphic to each other.) In particular, since
$k\in\{1,2,4,5,6\}$ we have $\dim H_1(\ee_k;\ZZ_p)\le1$ for any prime
$p\ne2,7$, and $\dim H_1(\ee_k;\ZZ_p)\le2$ for $p=2$ or $7$. As $X$ is
isomorphic to a quotient of $\ek$, it follows that $\dim
H_1(X;\ZZ_p)\le1$ for any prime $p\ne2,7$, and $\dim H_1(X;\ZZ_p)\le2$
for $p=2$ or $7$.  Hence by \fullref{refined SW corollary}, we
have $\dim H_1(M;\ZZ_p)\le2$ for $p\ne2,7$, and $\dim
H_1(M;\ZZ_p)\le3$ for $p=2,7$.
 
It remains to prove that if $\vol M\le\theotherbound$ then $\dim
H_1(M;\ZZ_7)\le2$. We have observed that in this case
$k\in\{2,4,5,6\}$. By the list of isomorphism types of the
$H_1(\ek;\ZZ)$ given above, it follows that $\dim
H_1(\ee_k;\ZZ_7)=0<1$. Hence in this case the argument given above for
$p\ne2,7$ goes through in exactly the same way to show that \mbox{$\dim
H_1(M;\ZZ_7)\le2$}.  \EndProof

\Proof[Proof of \fullref{main theorem}] For the case in which $M$
is non-exceptional, the theorem is an immediate consequence of
Propositions \ref{large R} and \ref{other main theorem}. For the case
in which $M$ is exceptional, the assertions of the theorem are
equivalent to those of \fullref{small R}.  \EndProof

\bibliographystyle{gtart}
\bibliography{link}

\end{document}